\documentclass[twocolumn,letterpaper]{IEEEtran}

\usepackage{cite} 
\usepackage{amsmath,amssymb,bm} 
\usepackage{bbm} 
\usepackage{graphicx} 
\usepackage[font=small]{caption} 
\usepackage[font=small]{subcaption} 
\usepackage{wrapfig} 
\usepackage{float} 
\usepackage{color,soul} 
\sethlcolor{red}
\usepackage{comment} 
\usepackage[hidelinks]{hyperref} 
\usepackage{cleveref} 
\crefname{subsection}{subsection}{subsections}
\Crefname{appsec}{appendix}{appendices}
\usepackage{tablefootnote} 
\usepackage{multirow} 
\usepackage{url} 
\usepackage{authblk} 
\usepackage{blindtext} 
\usepackage{footnote} 
\usepackage{caption}
\usepackage{subcaption}
\usepackage{algpseudocode}
\usepackage{courier}

\title{Near optimal tracking control of a class of nonlinear systems and an experimental comparison}

\author{Farshid Asadi\thanks{Department of Mechanical and Aerospace Engineering, Rutgers University, Piscataway, NJ 08854 USA (e-mail: fa416@scarletmail.rutgers.edu, farshidasadi47@yahoo.com).} and Ali Heydari\thanks{Mechanical Engineering, University of New Mexico, Albuquerque, NM 87131, USA.}}
\begin{document}

\maketitle

\hl{This article has been published in \textit{IET Control Theory \& Applications}. Content is final as presented, with the exception of pagination. Please cite the published version as,}

\ul{For Word}:

\noindent\texttt{Asadi, Farshid, and Ali Heydari. "Near optimal tracking control of a class of non‐linear systems and an experimental comparison." IET Control Theory \& Applications 14, no. 19 (2020): 3086-3096}

\ul{For Latex}:

\hspace{-18pt}\texttt{
@article\{asadi2020near,\\
  title=\{Near optimal tracking control of a class of non-linear systems and an experimental comparison\},\\
  author=\{Asadi, Farshid and Heydari, Ali\},\\
  journal=\{IET Control Theory \& Applications\},\\
  volume=\{14\},\\
  number=\{19\},\\
  pages=\{3086--3096\},\\
  year=\{2020\},\\
  publisher=\{Wiley Online Library\}\}
}
\vspace{10pt}

\begin{abstract}
In this paper, near optimal tracking of a class of nonlinear systems is addressed. Adaptive (approximate) dynamic programming approach is used to calculate the optimal control in closed form. ADP\footnote{Adaptive (approximate) dynamic programming} has been widely used to resolve optimal regulation and tracking problems of nonlinear control systems. Despite advances in the so called supervised and unsupervised ADP techniques for optimal tracking, they have a main draw back. That is, the optimal controller needs to be recalculated for every particular reference trajectory. The main goal of this work is to address this issue for a class of nonlinear systems. Finally, this approach is applied on a Delta robot and the performance of the method is analyzed experimentally.
\end{abstract}

\section{Introduction}
Trajectory tracking of nonlinear systems is a classic problem in control theory. Optimal control, as one of the approaches to solve this problem, has attracted some efforts throughout past several decades. The interested reader can refer to \cite{bryson1975applied,werbos2004adp,lewis2009reinforcement,saridis1979approximation,ccimen2008state,zhang2019near} for a concise introduction to more common nonlinear optimal control techniques. Since analytic solutions of nonlinear optimal control problems is not available, except for simple cases, seeking approximate solutions is a common practice. ADP technique \cite{werbos2004adp}, that is approximating cost function with a neural network and learning the optimal cost function in a backward manner (dynamic programming), is one of the widely used techniques among researchers. Optimal tracking problem have been studied for both continuous time and discrete time systems, but regardless of this, the pursued solutions can be categorized into two general frameworks: LQR\footnote{Linear quadratic regulator} extensions and ADP based approaches. 

In the first approach, the nonlinear plant is modeled as a linear system with time-varying matrices, and then techniques from linear optimal control theory are used. For instance, in \cite{lahdhiri1999design}, a feedback linearization is done on the nonlinear plant and then, a linear optimal problem is defined for the resultant feedback linearized system. In this approach, the object function is not directly related to the physical system and may not have physical realization. In \cite{ccimen2004nonlinear,jeyed2019development}, SDRE\footnote{State dependent Riccati equation} approach is used for input affine nonlinear systems. The main drawback of this method is that the proper choice of state-dependent quasilinear form plays an important role in the algorithm \cite{ccimen2008state}. Also in \cite{chou2004line}, a general nonlinear system is considered and error dynamics is estimated adaptively as a linear system. The optimal control then, is calculated based on the linear estimation.

In the ADP based approaches, the total cost is approximated with a function approximator of appropriate form and then, this approximation is used in order to calculate optimal cost and optimal control. ADP based approaches can be categorized into two branches based on the objective function that they used. In the first one, the objective function for tracking is defined based on the error and the total control input of the system. Optimizing this cost function leads to optimality but the resulted controller is not locally asymptotically stable in general, this will be discussed later. For instance, in \cite{tang2007optimal}, the general nonlinear system is decomposed based on its linearization and residual terms, then the optimal control is calculated as a combination of linear and residual part. In \cite{mclain1999synthesis}, a finite horizon continuous time optimal tracking is considered and then, the optimal control is calculated by direct implementation of ADP. A discrete time version of finite horizon approximate optimal tracking can also be found in \cite{heydari2014fixed}, in this work the controller can accept different initial conditions of the same reference trajectory dynamics. Also in \cite{modares2018adaptive,kiumarsi2014actor,modares2014optimal}, states and reference trajectory are augmented in a new variable and the optimal problem is solved as a regulation. In these three work reinforcement learning is used to calculate the optimal control online. Moreover in \cite{khiabani2019design}, a discrete time optimal tracking controller is considered for a switching system and is solved by direct implementation of ADP. Another ADP based reinforcement learning is used in \cite{yang2011reinforcement} for tracking control of a class of discrete time nonlinear systems with unknown dynamics, however in this approach it is assumed that the input transition matrix is positive definite. For problems with bounded input, \cite{lyshevski1999optimal} has proposed an approach by adding a non-quadratic functional to the total cost. This approach is not done with ADP for tracking problems. However, in \cite{abu2005nearly} it is used along with ADP for a regulation problem.

In the second category of ADP based approaches, the control input is decomposed into a steady state (that makes the error dynamic stationary at origin) and a transient part. Then, the control objective is defined based on the error and the transient control. This approach can be found in \cite{park1996optimal,zhang2008novel}. In \cite{zhang2008novel}, reference trajectory and error are augmented into new states and then the ADP is used to solve the resulted augmented regulation optimal problem. In \cite{dierks2009optimal}, an optimal control problem is defined for the transient control and then the total cost is calculated in an online manner by ADP. In this method knowledge of system dynamics is not necessary. In \cite{kamalapurkar2015approximate}, reinforcement learning is used to solve the optimal tracking control of a nonlinear system with unknown dynamics online. In these two works (\cite{dierks2009optimal,kamalapurkar2015approximate}) the optimal problem is solved for an augmented system, as in \cite{zhang2008novel}.

Despite advances in the mentioned works, ADP based methods all share a common drawback. That is, the controller needs to be re-calculated for each particular reference trajectory. This issue also exists in LQR based techniques.

This work is dedicated to solve the mentioned problem of ADP based approaches for a class of nonlinear systems. The proposed method uses the idea of control decomposition to eliminate trajectory dynamics from error dynamics, without eliminating systems dynamic matrices. The optimization is done based on using transient control in the objective function, which is called \textit{modified total cost}, in here. Effects of this decomposition on optimality and asymptotic stability of the closed loop system will be discussed, which, to the best of our knowledge, has not been done in the related literature yet. Furthermore, it will be shown that by optimizing expectation of modified total cost, instead of its exact value, there is no need to know the reference trajectory in the training stage. This change will lead to the main contribution of this paper, that is a near optimal asymptotically stabilizing tracking controller that does general tracking for a class of nonlinear systems. Finally, it will be shown that using optimal control based on expected value of modified total cost (instead of its exact value), does not hurt asymptotic stability of the closed loop system. The proposed controller is near optimal in three aspects. First, because of the form of steady state control that is used. Second, it optimizes expected value of  modified total cost, instead of exact modified total cost of a reference trajectory. Third, it approximates this objective function which is the core of ADP.

In what follows, first the problem is defined and the proposed method is explained. Then theoretical support for optimality, convergence, and asymptotic stability of the approach is presented. Finally, the method is implemented experimentally on a Delta parallel manipulator and its performance is shown in comparison to some standard nonlinear control techniques.
\section{Problem statement and resolution}
Let us define a tracking problem for a nonlinear system of the following form
\begin{equation}
x^{(p)} = f(X) + g(X)u,
\label{eq:csys}
\end{equation}
where $x(t) \in \mathbb{R}^{m\times 1}$ is an $m\mbox{-}$vector\footnote{All vectors are column vectors.} of output of interest, $X = [x^\intercal, \cdots, {x^{(p-1)}}^\intercal]^\intercal \in \mathbb{R}^{(n\times 1)}$ is $n\mbox{-}$vector of states, and $u(t) \in \mathbb{R}^{m\times 1}$ is $m\mbox{-}$vector of control input. Also note that $n = m\times p$. Furthermore, $f(.): \mathbb{R}^{n\times1}\rightarrow \mathbb{R}^{m\times1}$ and $g(.): \mathbb{R}^{n\times1}\rightarrow \mathbb{R}^{m\times m}$ are functions representing dynamics of the system. Moreover $f(.)$ and $g(.)$ and their Jacobians are assumed to be continuous. Also, $x^{(i)}$ denotes the $i_{th}$ time derivative of $x$. The system is supposed to follow a particular reference trajectory (not known a priori), that is $X_d(t) = [x_d^\intercal, \cdots, {x_d^{(p-1)}}^\intercal]^\intercal \in D\subset \mathbb{R}^{n\times 1}$, with zero tracking error. Furthermore, assume that $\dot{X}_d$ exists and is continuous. Tracking error is defined as $E = [e^\intercal, \cdots, {e^{(p-1)}}^\intercal]^\intercal$ where $e = x-x_d$. As mentioned in the introduction, in some of the related literature, the optimal controller is designed to minimize the following cost function
\begin{equation}
J =  \int_{0}^{\infty}exp(-\rho \tau)(E^\intercal QE+{u}^\intercal R{u})d\tau,
\label{eq:ctotcost}
\end{equation}
where $Q \in \mathbb{R}^{n\times n},R \in \mathbb{R}^{m\times m}$, and $\rho \in \mathbb R_{\ge 0}$ are semi-positive definite error penalizing matrix, positive definite control penalizing matrix, and discount factor, respectively. The error dynamics of this tracking problem can be written in the following form
\begin{equation}\label{eqn_cedyn}
e^{(p)} = f(E+X_d)+g(E+X_d)u-x_d^{(p)}.
\end{equation}
This error dynamics can be also stated in state space from as
\begin{equation}\label{eqn_csedyn}
\dot{E} = F(E+X_d)+G(E+X_d)u-\dot{X}_d,
\end{equation}
where one has
\begin{equation*}
F(E+X_d) = \begin{bmatrix}
\dot{e} + \dot{x}_d\\
\vdots\\
e^{(p-1)} + x_d^{(p-1)}\\
\\ f(E+X_d)
\end{bmatrix},
\end{equation*}
\begin{equation*}
G(E+X_d) = \begin{bmatrix}
0_{(n-m) \times m}\\ g(E+X_d)
\end{bmatrix}.
\end{equation*}
Note that the above dynamics is non-autonomous\footnote{The dynamic system has a direct dependence on time through $X_d$ and its time derivative.} and non-stationary\footnote{The equilibrium point of the dynamic system does not lie at origin.} at origin with respect to its states $E$.

Even though error dynamics and objective function in the form of \cref{eqn_cedyn,eq:ctotcost} are commonly used, there are three disadvantages with formulating the problem in this way. First, the optimal tracking controller is not locally asymptotically stabilizing in general (see \cref{appendixexample}). The reason is that the reference trajectory, generally is not an invariant set of the system dynamics, which is needed for optimal control to be asymptotically stabilizing. This shows itself as a steady state error\footnote{This steady state error is because of the analytical construction of the controller, not from disturbance and\textbackslash or uncertainties.}. Second, the presence of discounting factor $\rho$ means that, just a limited part of horizon is important to the controller. This will lead to a higher steady state error and worsens the effects of the first problem. Third, the resulted optimal control can only follow\footnote{Assuming that the application tolerates the first and the second mentioned problems.} the reference trajectory that is solved for. This means that for new reference trajectories, the problem should be re-solved. These issues have motivated some authors to use a modified objective function and to revisit components of the error dynamics by decomposing the control input into \textit{steady} and \textit{transient} parts. However, the third issue is not solved in any of the ADP related literature yet, to the best of the authors' knowledge. Furthermore, interpretation of such control decomposition with respect to optimality is not done in the referenced works, to the best of our knowledge.

In a tracking problem, the evolution of the system can be categorized in two phases: the transient and the steady state. This gives an idea of decomposing the control to a steady state control plus a correction term, when it is possible. For a system in the form of \cref{eq:csys}, the steady state control, that is $u_s$, can be defined to satisfy the following equation
\begin{equation}\label{eqn_cpresscont}
x_d^{(p)} = f(X_d)+g(E+X_d)u_s,
\end{equation}
this form of steady state control is used in \cite{dierks2009optimal} for a discrete-time system, and its main advantage over other forms in literature is that it eliminates trajectory dynamics from error dynamics of \cref{eqn_cedyn}. A controllable plant is assumed, therefore $g(E+X_d) = g(X)$ is invertible. Then $u_s$ can be calculated as
\begin{equation}\label{eqn_csscont}
u_s = g^{-1}(E+X_d)(x_d^{(p)}-f(X_d)).
\end{equation}
At any instance, if error equals zero, that is $E=0$, then applying $u_s$ leads to perfect tracking. The total control is the sum of steady state control, that is $u_s$, and a corrective term, that is $\Delta u$, so it is defined as
\begin{equation}
u = u_s + \Delta u.
\label{eq:ctotcont}
\end{equation}
By substituting \cref{eq:ctotcont,eqn_csscont} in \cref{eqn_cedyn}, the error dynamics equation, (\cref{eqn_cedyn}) can be rewritten as
\begin{equation}\label{eqn_cedynred}
e^{(p)} = f(E+X_d)-f(X_d)+g(E+X_d)\Delta u,
\end{equation}
furthermore, this equation can be written in state space form as
\begin{equation}\label{eqn_csedynred}
\dot{E} = F(E+X_d)-F(X_d)+G(E+X_d)\Delta u.
\end{equation}
The above error dynamics is stationary at origin. Now, one can define modified total cost based on the corrective term $\Delta u$ (instead of total control $u$) as
\begin{equation}
V = \int_{0}^{\infty}(E^\intercal QE+{\Delta u}^\intercal R{\Delta u})d\tau,
\label{eq:ctotcostred}
\end{equation}
where optimal transient control will be calculated by optimizing the above total cost. This cost function is commonly used in this category of solutions to the optimal tracking problem. Decomposing the control and redefining the total cost resolves the first mentioned problem. Assuming that the system is controllable, the above modified total cost is bounded. The reason is that $\Delta u$ vanishes as the transient phase finishes. This cost function minimizes the error and the corrective control term, so it brings the system to the steady state tracking phase asymptotically. The reason is that, by imposing the steady state control and optimizing modified total cost, the optimal problem is converted to an optimal regulation problem which is asymptotically stable (see for example \cite{lyashevskiy1995control} for asymptotic stability of optimal regulation problem). Furthermore, since the boundedness of modified total cost, that is \cref{eq:ctotcostred}, is achieved without introducing a discounting factor, there is no risk of associated steady state error. Therefore, the second issue is also resolved.

For any particular reference trajectory in time, modified total cost, that is \cref{eq:ctotcostred}, only depends on the initial error, that is $E_0$. This is a key point in this analysis that also reduces dimensionality of the value function and therefore, mitigates curse of dimensionality further. However, the issue is that the trajectories are not known ahead of time. If one writes HJB\footnote{Hamilton- Jacobi- Bellman} equation for \cref{eq:ctotcostred} it can be seen (see \cref{Theoretical background}) that knowledge of trajectory is needed for calculating the modified optimal cost function. One solution to this issue is using expectation of total cost instead of its exact value for a specific trajectory in time. The reason is that one can consider the desired trajectory, that is $X_d(t)$, as a parameter with uniform distribution in ROI\footnote{Region of interest}. This leads to the main contribution of this paper. The expected value of modified total cost, that is \cref{eq:ctotcostred}, can be written in the following form,
\begin{equation}
\overline{V}(E_0) = \underset{X_d\in ROI}{\mathbb{E}}
\left\{\int_{0}^{\infty}(E^\intercal QE+{\Delta u}^\intercal R{\Delta u})d\tau)\right\},
\label{eq:cstochtotcostred}
\end{equation}
where $\mathbb{E}$ denotes mathematical expected value. For every specific trajectory and the defined problem, optimal modified total cost exists uniquely \cite[pp. 284--291]{athans2013optimal} and is two times differentiable \cite{strulovici2015smoothness}. Therefore, its expected value (that is simply an average over all possible trajectories in the present case) also exists, is unique and is two times differentiable. By following the procedure of \cite[pp. 131--136]{bryson1975applied} and taking time derivative of \cref{eq:cstochtotcostred} and using the error dynamics from \cref{eqn_csedynred}, the non-optimal HJB equation can be derived as the following 
\begin{equation}
\underset{X_d\in ROI}{\mathbb{E}}\left\{\overline{V}^\intercal_E(F-F_d+G\Delta u)+E^\intercal QE+\Delta u^{\intercal}R \Delta u \right\}= 0,
\label{eq:ctotcostrecredhjb}
\end{equation}
where $\overline{V}_E = \frac{d\overline{V}(E)}{dE}$, $F = F(E+X_d) = F(X)$, $F_d = F(X_d)$, and $G = G(E+X_d) = G(X)$.
Optimal transient control, that is $\Delta u^*$, is the minimizer of LHS\footnote{Left hand side} of \cref{eq:ctotcostrecredhjb}, and can be calculated as
\begin{equation}
{\Delta u}^* = -\frac{1}{2}R^{-1} G^\intercal \overline{V}_E^*,
\label{eq:coptcontred}
\end{equation}
where $\overline{V}_E^*(E)$ is gradient of expectation of optimal modified total cost.

There are several ways in the literature to solve the resulted HJB equation, including PI\footnote{Policy iteration} algorithm \cite{leake1967construction,saridis1979approximation}, integral PI \cite{beard1997galerkin}, integral VI\footnote{Value iteration} \cite{bian2016value}, projection technique \cite{kompas2010comparison}, perturbation method \cite{kompas2010comparison}, and parametric linear programming technique \cite{kompas2010comparison}. Among these methods integral VI is chosen. The reason is that it gives a good understanding of underlying theory and it does not need initial admissible policy as needed in other iterative methods based on PI. One can rewrite \cref{eq:cstochtotcostred} in the following form and apply Bellman principle of optimality \cite{Bellman:1957} as
\begin{multline}
\overline{V}^*(E(t)) =\underset{X_d\in ROI}{\mathbb{E}}
\left\{\int_{t}^{t+\Delta T}(E^\intercal QE+{\Delta u^*}^\intercal R{\Delta u^*})d\tau)
\right.\\ \left.+ \overline{V}^*(E(t+\Delta T))\right\},
\label{eq:cstochtotcostredrec}
\end{multline}
where $\Delta T \rightarrow 0$. This equation can be used in the so called integral value iteration to learn the expectation of modified optimal cost (in the other words, value function), from the following iterative procedure,
\begin{multline}
\overline{V}_{i+1}(E(t)) =\\ \underset{X_d\in ROI}{\mathbb{E}}
\left\{\int_{t}^{t+\Delta T}(E^\intercal QE+{\Delta u_{i+1}^*}^\intercal R{\Delta u_{i+1}^*})d\tau)
\right.\\\left.+ \overline{V}_i(E(t+\Delta T))\right\},
\label{eq:cstochvfredrec}
\end{multline}
where $\Delta u_{i+1}^*$ is defined as
\begin{equation}
{\Delta u}_{i+1}^* = -\frac{1}{2}R^{-1} G^\intercal \overline{V}^*_{i_E}.
\label{eq:coptcontredrec}
\end{equation}
Note that because of using the expectation of total cost, there is no need for knowing the trajectory in training stage. This means that once the expectation of modified value function is calculated, it can be used to track every trajectory in the ROI. Therefore the third issue with existing ADP based methods is also solved for nonlinear systems with the dynamics given by \cref{eq:csys}.

To calculate expected modified value function, and consequently optimal transient control, in a closed form, ADP \cite{werbos2004adp} is used here. To do this, value functions in \cref{eq:cstochtotcostredrec,eq:cstochvfredrec} are approximated with a linear (in weight) NN\footnote{Neural network} as
\begin{equation}
\overline{V}(E) = W^\intercal \varphi (E),
\label{eq:cvfapp}
\end{equation}
\begin{equation}
\overline{V}_i(E) = W_i^\intercal \varphi (E),
\end{equation}
where $\varphi (E)$ is the vector of basis functions.

The procedure for training the neural network with integral value iteration can be summarized as:
\begin{enumerate}
	\item[$a$.] Initialize some random values for $X$ and $X_d$ in the ROI.
	\item[$b$.] Calculate errors $E$, based on values generated in step $a$.
	\item[$c$.] Initialize \cref{eq:cstochvfredrec} with $\overline{V}_{0}(E) = 0$.
	\item[$d$.] For every $E$ and $X_d$ pair and some small constant value of $T$, repeat the following stages:
	\begin{enumerate}
		\item[1.] Calculate ${\Delta u}_{i+1}^*(E,X_d)$ from \cref{eq:coptcontred}.
		\item[2.] Calculate $\overline{V}_{i+1}(E)$ from \cref{eq:cstochvfredrec}.
	\end{enumerate}
	\item[$e$.] Use values from step $d$ to update weights from the least square method \cite[pp. 302--305]{boyd2004convex}.
	\item[$f$.] Calculate $\delta = max(\left|W_{i+1}-W_i\right|)$ and do one of the following stages:
	\begin{enumerate}
		\item[1.] If $\delta< threshold$, terminate the iterative procedure and set $W = W_{i+1}$ and go to step g.
		\item[2.] If $\delta \ge threshold$ then go to step d.
	\end{enumerate}
	\item[$g$.] The optimal corrective term can be calculated as $${\Delta u}^* = -\frac{1}{2}R^{-1} G^\intercal(E+X_d) W^\intercal \Delta \varphi(E),$$ where $\Delta \varphi(E) = \frac{\partial \varphi(E)}{\partial E}$ is gradient of $\varphi(E)$.
\end{enumerate}
\section{Theoretical background}\label{Theoretical background}
\subsection{Optimality}\label{sec:optimality}
In the presented approach, the total control is calculated as a combination of the steady state control and the optimal transient control. Considering this decomposition, one may ask about the optimality of the resulted total cost with respect to the original cost function, that is \cref{eq:ctotcost}. To make the analysis easier to grasp, we will investigate the case of a specific reference trajectory to avoid getting involved in the expected values of the total costs. Let us define a new steady state control $u_r$ and a discounted version of \cref{eq:ctotcostred} as
\begin{equation}\label{eqn_crrcont}
u_r = g^{-1}(X_d)(x_d^{(p)}-f(X_d)),
\end{equation}
\begin{equation}
W(E) = \int_{0}^{\infty}exp(-\rho \tau)(E^\intercal QE+{\Delta u}^\intercal R{\Delta u})d\tau.
\label{eq:cstochtotcostreddisc}
\end{equation}
By substituting total control based on $u_r$, that is $u = u_r + \Delta u$, in \cref{eq:ctotcost}, the total cost can be rewritten as
\begin{multline}
J(E) = \int_{0}^{\infty}exp(-\rho \tau)(E^\intercal QE + \Delta u^\intercal R \Delta u \\+ 2{\Delta u}^\intercal R{u_r} + u_r^\intercal R u_r)d\tau,
\label{eq:ctotcostn1}
\end{multline}
or
\begin{equation}
J(E) = W(E) + \int_{0}^{\infty}exp(-\rho \tau)(2{\Delta u}^\intercal R{u_r} + u_r^\intercal R u_r)d\tau.
\label{eq:ctotcostn2}
\end{equation}
By minimizing \cref{eq:ctotcostn1,eq:ctotcostn2} over $\Delta u$ and eliminating $J^*(E)$ among them one can write
\begin{multline}
\min_{\Delta u}\left\{W(E) + \int_{0}^{\infty}exp(-\rho \tau)(2{\Delta u}^\intercal R{u_r} + u_r^\intercal R u_r)d\tau\right\}\\=
\min_{\Delta u}\biggl\{\int_{0}^{\infty}exp(-\rho \tau)(E^\intercal QE + \Delta u^\intercal R \Delta u \\+ 2{\Delta u}^\intercal R{u_r} + u_r^\intercal R u_r)d\tau\biggr\}.
\label{eq:ctotcostn3}
\end{multline}
The minimization is on the same variable on both sides of \cref{eq:ctotcostn3}, so it can be simplified as
\begin{equation}
\begin{split}
\min_{\Delta u}\left\{\right.&\left.W(E)\right\} = W^*(E) =\\& \min_{\Delta u}\left\{\int_{0}^{\infty}exp(-\rho \tau)(E^\intercal QE + \Delta u^\intercal R \Delta u)d\tau\right\}.
\end{split}
\label{eq:ctotcostn4}
\end{equation}
The process above means that optimizing $W(E)$ is equal to optimizing $J(E)$, while imposing $u_r$ on the system. Furthermore, as $\rho \rightarrow 0$, $W(E)\rightarrow V(E)$. This means that optimizing $V(E)$ is equivalent of optimizing $W(E)$ while imposing $u_r$ on the system and $\rho \rightarrow 0$. Since imposing $u_r$ on the error dynamics and finding optimal transient control, for a specific reference trajectory, transforms the optimal control problem to an equivalent optimal regulation problem, so it is asymptotically stable (see \cite{lyashevskiy1995control} for asymptotic stability of optimal regulation problems). Also since $u_r$ is the exact steady state control, the resulted total control is optimal among asymptotically stabilizing controllers, however it is not the absolute optimal control (that as discussed, see \cref{appendixexample}, is not generally asymptotically stabilizing). Furthermore, because $u_s$ (that is used in the proposed method) has a slight difference with $u_r$, the proposed method has a degree of sub-optimality. Therefore, the proposed method, for a specific reference trajectory, is a near optimal control among asymptotically stabilizing tracking controllers.

The other point that should be investigated is about optimality of optimal transient control, that is \cref{eq:coptcontred}, with respect to expected value of modified total cost, that is \cref{eq:cstochtotcostred}. In the presented approach, the optimization problem is defined based on the expectation of modified total cost. Introducing the expectation in the equations, a legitimate question is the optimality of the selected control, that is \cref{eq:coptcontred}. To answer this, one needs to take derivative of LHS of \cref{eq:ctotcostrecredhjb} with respect to $\Delta u$ and solve the following equation for $\Delta u$
\begin{multline}
\frac{\partial}{\partial \Delta u}\biggl\{\underset{X_d\in ROI}{\mathbb{E}}\Bigl\{\overline{V}^\intercal_E(F-F_d+G\Delta u)\\+E^\intercal QE+\Delta u^{\intercal} R\Delta u \Bigr\}\biggr\}  = 0.
\label{eq:ctotcostrechjbdif}
\end{multline}
In the above equation, since the expectation is over $X_d$, the derivative can be interchanged with the expectation and one has
\begin{multline}
\underset{X_d\in ROI}{\mathbb{E}} \biggl\{\frac{\partial}{\partial \Delta u}\Bigl\{\overline{V}^\intercal_E(F-F_d+G\Delta u)\\+E^\intercal QE+\Delta u^{\intercal} R\Delta u \Bigr\}\biggr\}  = 0,
\label{eq:ctotcostrechjbdifin}
\end{multline}
note that if the following holds
\begin{equation}
\begin{split}
\frac{\partial}{\partial \Delta u}\left\{\overline{V}^\intercal_E(F-F_d+G\Delta u)+E^\intercal QE+\Delta u^{\intercal} R\Delta u \right\}  = 0,
\end{split}
\label{eq:ctotcostrechjbdifinred}
\end{equation}
then also \cref{eq:ctotcostrechjbdifin} holds, therefore the answer to the above equation, which is well known to be in the form of \cref{eq:coptcontred}, is a solution to \cref{eq:ctotcostrechjbdif}. Consequently, \cref{eq:coptcontred} is a minimizer to \cref{eq:cstochtotcostred,eq:ctotcostrecredhjb}. Therefore, optimality of \cref{eq:coptcontred} with respect to the expectation of modified total cost is proved.
\subsection{Convergence}
Since integral value iteration of \cref{eq:cstochvfredrec,eq:coptcontredrec} is an iterative procedure, convergence of the iterations is a concern. This concern is investigated in this subsection while neglecting approximation error of \cref{eq:cvfapp}. The procedure is adopted from \cite{bian2016value} with some changes. Let us define the following transformations, for some positive $V(E)$ with $V(0) = 0$, as
\begin{multline}
T^{\Delta u}(V,X_d) =
\int_{t}^{t+\Delta T}(E^\intercal QE+{\Delta u}^\intercal R{\Delta u})d\tau)
\\+ V(E(t+\Delta T)),
\label{eq:transdef}
\end{multline}
\begin{equation}
\overline{T}^{\Delta u}(V) = \underset{X_d\in ROI}{\mathbb{E}}\left\{T^{\Delta u}(V,X_d)\right\}.
\label{eq:etransdef}
\end{equation}
Also let $T(V,X_d)$ be defined as $T^{^{\Delta u^{\circ}}}(V,X_d)$, that is calculated with ${\Delta u}^{\circ} = -\frac{1}{2}R^-1G(E+X_d)^\intercal V_E$ which is the minimizer of right hand side of \cref{eq:transdef} (and also \cref{eq:etransdef}, in the same way explained in \cref{sec:optimality}). Furthermore, let $\overline{T}(V)$ be the expected value of $T(V,X_d)$. In this way one can write \cref{eq:cstochtotcostredrec,eq:cstochvfredrec} in the forms of $\overline{V}^* = \overline{T}(\overline{V}^*)$ and $\overline{V}_{i+1} = \overline{T}(\overline{V}_i)$, respectively.

If for some $V_l$ and $V_k$ the inequality $V_l< V_k$ holds, then $T^{\Delta u}(V_l,X_d)< T^{\Delta u}(V_k,X_d)$ also holds for all $E, X_d \in ROI \text{ and } E\neq 0$. Therefore one can conclude that $\overline{T}^{\Delta u}(V_l)< \overline{T}^{\Delta u}(V_k)$. Consequently, by this assumption, one has $\overline{T}(V_l)\leq \overline{T}^{{\Delta u}^{\circ}_k}(V_l)< \overline{T}(V_k)$).

If $\overline{V}_0 = 0$ holds, by investigating \cref{eq:cstochvfredrec} one can see that $\overline{V}_1> \overline{V}_0$. Also if one assumes $\overline{V}_{i+1}>\overline{V}_i$, then $\overline{V}_{i+2} = \overline{T}(\overline{V}_{i+1})>\overline{V}_{i+1} = \overline{T}(\overline{V}_i)$. Thus, if integral value iteration starts with $\overline{V}_0 = 0$, then by induction one has $\overline{V}_{i+1}>\overline{V}_i$ for all $E, X_d \in ROI \text{ and } E\neq 0$.

Since a controllable plant is assumed, then there exist a  non-optimal (in the sense of expectation of modified total cost) stabilizing control policy, that is $\Delta h$, whose expected modified total cost, that is $\overline{Z}_0(e)$, is greater than $\overline{V}_0$ and $\overline{V}^*$. Also note that one can write $\overline{Z}_0 = \overline{T}^{\Delta h}(\overline{Z}_0)$. Therefore $\overline{Z}_1 = \overline{T}(\overline{Z}_0)< \overline{Z}_0 = \overline{T}^{\Delta h}(\overline{Z}_0)$. Consequently one can conclude $\overline{Z}_{i+1}< \overline{Z}_i$, in the same manner used in the previous paragraph.

Furthermore, since $\overline{V}^*<\overline{Z}_0$, then $\overline{V}^* = \overline{T}(\overline{V}^*)< \overline{Z}_1 = \overline{T}(\overline{Z}_0)$. By repeating this for $i$ times, one can conclude that $V^* < Z_i$ . Since $Z_i$ is a decreasing positive sequence\footnote{Elements of the sequence are positive.} that is lower bounded by $V^*$, it converges to this lower-bound. Moreover, since $\overline{V}_0<\overline{Z}_0$, by the same reasoning $\overline{V}_i<\overline{Z}_i$.

Finally, assuming that $\overline{V}_0 = 0$, one can combine strictly monotonic convergence of $\overline{Z}_i$ to $\overline{V}^*$, strictly monotonic increase of $\overline{V}_i$, the inequality $\overline{V}_i<\overline{Z}_i$, and positiveness of these functions to conclude that $\overline{V}_i$ also converges to $\overline{V}^*$ as $i \rightarrow \infty$ for all $E, X_d \in ROI \text{ and } E\neq 0$. Furthermore, one has $\overline{V}_i(0) = \overline{Z}_i(0) = \overline{V}^*(0) = 0$ from their construction. As a result, convergence of the integral value iteration for the proposed method is guaranteed for all $E, X_d \in ROI$.
\subsection{Stability}
The most important aspect of a controller is its stability. We will show asymptotic stability of the proposed method through an appropriate Lyapunov function. Let us define the \textit{Hamiltonian}, for any differentiable function $V_l(E): \mathbb{R}^{n\times 1}\rightarrow \mathbb{R}$, as
\begin{equation}
H_{\Delta u}(V_l) := V_{l_E}^\intercal(F-F_r + G\Delta u) + E^\intercal QE + {\Delta u}^\intercal R\Delta u.
\label{eq:hamiltonian}
\end{equation}
Furthermore, let $H(V_l)$ be optimal value of $H_{\Delta u}(V_l)$ with respect to $\Delta u$. Moreover, let us define $Z(E)$ as the modified total cost of an admissible control\footnote{In this work, admissible controls are limited to asymptotically stabilizing controllers.}, calculated from \cref{eq:ctotcostred} for a specific reference trajectory. Based on Lemma 1 of \cite{leake1967construction}, if the following inequality holds, for any $V_l$,
\begin{equation}
H(Z) \leq H(V_l),
\label{eq:hamil1}
\end{equation}
then one has
\begin{equation}
V_l(E) \leq Z(E).
\label{eq:hamil2}
\end{equation}
One can write optimal value of the modified total cost of a specific reference trajectory, that is calculated form \cref{eq:ctotcostred} (here called $S^*(E)$ for clarity of notation), as
\begin{equation}
S^*(E) = \overline{V}^*(E) + D^*(E).
\label{eq:totdec}
\end{equation}
Furthermore, since $S^*$ is the solution of HJB equation, so one has
\begin{equation}
H(S^*) = {S^*_E}^\intercal(F-F_r + G\Delta u^*_{s}) + E^\intercal QE + {\Delta u^*_{s}}^\intercal R\Delta u^*_{s} = 0,
\label{eq:hamilopt}
\end{equation}
where $\Delta u^*_{s} = \Delta u^* + \Delta u^*_{d} = -\frac{1}{2}R^{-1}G^\intercal (\bar{V}_E^* + D^*_E)$. This equation can be rewritten as
\begin{multline}
{D^*_E}^\intercal(F-F_r + G(\Delta u^* + \Delta u^*_{d}))\\ + E^\intercal QE + (\Delta u^* + \Delta u^*_{d})^\intercal R(\Delta u^* + \Delta u^*_{d})\\  + {\overline{V}^*_E}^\intercal(F-F_r + G(\Delta u^* + \Delta u^*_{d}))= 0.
\label{eq:rehamilopt}
\end{multline}
For any specific reference trajectory, there could be three cases, based on $D^*(E)$. First, assume that $D^* = 0$. In this cases $\overline{V}^*(E) = V^*(E)$. So, by substituting \cref{eq:totdec} into \cref{eq:hamilopt} one has
\begin{equation}
\dot{\overline{V}}^* = {\overline{V}_E^*}^\intercal(F-F_r + G\Delta u^*) = - E^\intercal QE - {\Delta u^*}^\intercal R\Delta u^*.
\label{eq:dl1}
\end{equation}
RHS of the above equation is negative definite (from definition of the problem). Therefore, when $D^* = 0$ and the system is controlled by the proposed controller, that is $\Delta u^*$ from  \cref{eq:coptcontred}, $\dot{\overline{V}}^*$ is negative definite.

Second, assume that $D^*(E)<0$ for $E\neq0$. In this case, negative definiteness of the $\dot{\overline{V}}^*$ under the proposed controller can be proved by contradiction. Assume that $\dot{\overline{V}}^*$ is negative except for some $E, X_d \in ROI \textit{ and } E\neq0$. Therefore one can write $\dot{\overline{V}}^* = {\overline{V}_E^*}^\intercal(F-F_r + G\Delta u^*)\geq0$ for some $E, X_d \in ROI \textit{ and } E\neq0$, which leads to $H(\overline{V}^*)\geq H(V^*) = 0$. In this way one can conclude that $\overline{V}^*\leq V^*$ (from Lemma 1 of \cite{leake1967construction}) for some $E, X_d \in ROI \textit{ and } E\neq0$. This is contradictory with the assumption of $D^*(E)<0$, so one cannot have $\dot{\overline{V}}^*\geq 0$ for some $E, X_d \in ROI \textit{ and } E\neq0$. From the same reasoning, one also cannot have $\dot{\overline{V}}^*\geq 0$ for all $E, X_d \in ROI \textit{ and } E\neq0$. Therefore, $\dot{\overline{V}}^*$ is negative definite under the proposed controller for $D^*(E)<0$, taking into account that $\dot{\overline{V}}^* = 0$ for $E = 0$ from its construction.

Third, assume that $D^*(E)>0$ for $E\neq0$. In this way, One has $H(D^*)> H(S^*) = 0$ for all $E, X_d \in ROI , E\neq0$. The reason is that assuming $H(D^*)\leq H(S^*)$ holds for some $E, X_d \in ROI \textit{ and } E\neq0$, leads to $D^*\geq S^*$, which is contradictory to $S^* = \overline{V}^* + D^*$ in the current case. Also one has $H(D^*) = 0$ for all $X_d \in ROI , E=0$, by its definition.

Note that one can write
\begin{equation}
{D^*_E}^\intercal G\Delta u^* = {\overline{V}^*_E}^\intercal G\Delta u^*_d = -2{\Delta u^*}^\intercal R\Delta u^*_d.
\label{eq:aux1}
\end{equation}
By expanding \cref{eq:rehamilopt}, one can write
\begin{multline}
{D^*_E}^\intercal(F-F_r + G\Delta u^*_{d})+ {D^*_E}^\intercal G\Delta u^* + E^\intercal QE \\+ {\Delta u^*}^\intercal R\Delta u^* + 2{\Delta u^*}^\intercal R\Delta u^*_d +{\Delta u^*_{d}}^\intercal R\Delta u^*_{d}\\  + {\overline{V}^*_E}^\intercal G\Delta u^*_d + {\overline{V}^*_E}^\intercal(F-F_r + G\Delta u^*)= 0.
\label{eq:rehamiloptexpand}
\end{multline}
By substituting \cref{eq:aux1} into \cref{eq:rehamiloptexpand} and combining the terms one has
\begin{multline}
{D^*_E}^\intercal(F-F_r + G\Delta u^*_{d}) + E^\intercal QE \\+ (\Delta u^* - \Delta u^*_{d})^\intercal R(\Delta u^* - \Delta u^*_{d})\\+ {\overline{V}^*_E}^\intercal(F-F_r + G\Delta u^*)= 0,
\label{eq:aux2}
\end{multline}
furthermore, by minimizing first three terms in \cref{eq:aux2}, it can be written as
\begin{multline}
{D^*_E}^\intercal(F-F_r + G\Delta u^*_{d}) + E^\intercal QE \\+{\Delta u^*_{d}}^\intercal R\Delta u^*_{d}+ {\overline{V}^*_E}^\intercal(F-F_r + G\Delta u^*)\leq 0.
\label{eq:aux3}
\end{multline}
First three terms of \cref{eq:aux3} are equal $H(D^*)$, by definition. By using $H(D^*)>0$ for $E\neq0$ (as per current case) in \cref{eq:aux3}, one will conclude that for all $X_d,E\in ROI, E\neq0$,
\begin{equation}
{\overline{V}^*_E}^\intercal(F-F_r + G\Delta u^*)< 0.
\label{eq:aux4}
\end{equation}
Also one has ${\overline{V}^*_E}^\intercal(F-F_r + G\Delta u^*) = 0$ for $E=0$, by definition. As a result, $\dot{\overline{V}}^*$ is negative definite under the control $\Delta u^*$ for $D^*(E)>0$.

$\overline{V}^*(E)$ is positive definite by its construction and as we proved $\dot{\overline{V}}^* = {\overline{V}_E^*}^\intercal(F-F_r + G\Delta u^*)$ is negative definite in the ROI. Therefore $\overline{V}^*(E)$ is a Lyapunov function for the system under control of $\Delta u^*$, that is from \cref{eq:coptcontred}. Therefore, the closed loop system, from the proposed method, is locally asymptotically stable.
\section{An experimental case study}
\begin{figure}[H]
	\centering
	\includegraphics[width=.48\textwidth]{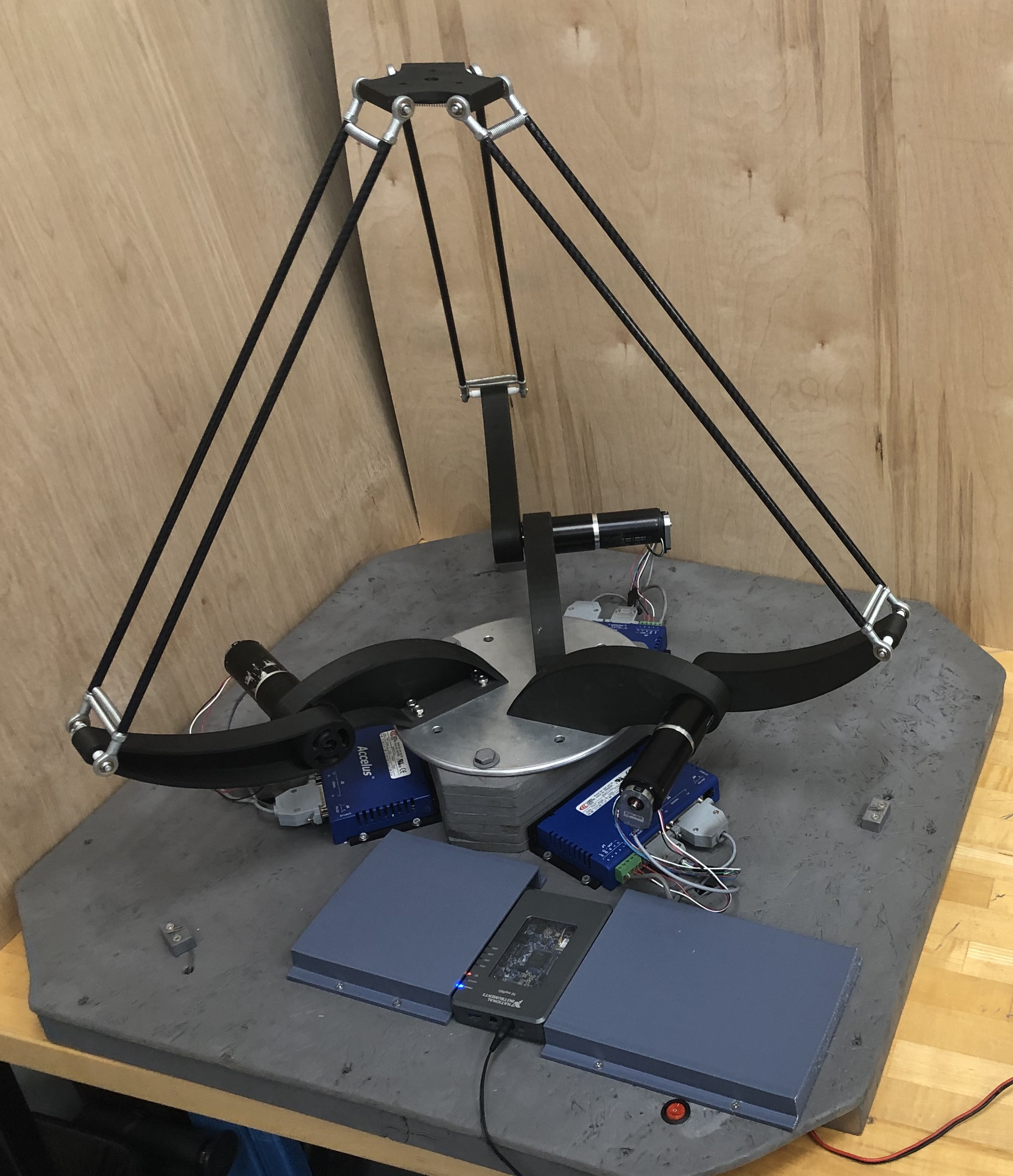}
	\caption{Delta robot used in the experiments}\label{fig:delta}
\end{figure}

To show and compare the performance of the presented approach, an experimental study is done on a developed Delta parallel robot, as depicted in \cref{fig:delta}.

Delta robot is a parallel manipulator with three transnational DOFs\footnote{Degrees of freedom} designed by Clavel \cite{clavelphdthesis}. The dynamic model of this robot can be presented in the following form \cite{tsai1999robot}
\begin{equation}
M(w,q)\ddot{w} + C(w,q,\dot{w},\dot{q})\dot{w} + G(w,q) = J(w,q)^\intercal \tau,
\label{eq:deldyn}
\end{equation}
where $M,C,G,w,q$, and $\tau$ are mass matrix, Coriolis matrix, gravitational vector, workspace coordinate vector, joint space coordinate vector, and motor torques, respectively. Different methods have been used to control the Delta robot \cite{castaneda2014robust,paccot2009review,taghirad2013parallel,codourey1998dynamic}. In the present work, computed torque method,which is usually used as a benchmark in the related literature, and sliding mode control, which is suitable benchmark of robustness, are considered for comparison.

The control law based on CT\footnote{Computed torque} can be written as \cite{taghirad2013parallel}
\begin{equation}
\tau = J^{-\intercal}\left(C\dot{w} + G + M(\ddot{w}_d-K_d\dot{\tilde{w}}-K_p\tilde{w} )\right),
\label{eq:cct}
\end{equation}
where $w_d,K_d$, and $K_p$ are desired position, derivative gain, and proportional gain, respectively. $\tilde{w} = w-w_d$ is also position error.

The control law for SMC\footnote{Sliding mode control} can also be written as \cite{slotine1991applied}
\begin{equation}
\tau = J^{-\intercal}\left(C\dot{w} + G + M(\ddot{w}_d-\lambda \dot{\tilde{w}}-Ksat(\frac{S}{\varphi}) )\right),
\label{eq:csmc}
\end{equation}
where $S = \dot{\tilde{w}} + \lambda \tilde{w}$, $\lambda$, $K$, and $\varphi$ are sliding surface, sliding surface parameter, sliding mode controller gain, and boundary layer, respectively. Also $sat$ denotes saturation function.

To train ADP based controller, $500$ sets of randomly generated data is used and the training is done $10$ times independently with least square method. Then weights of these 10 trainings are averaged and used for the experiments. Based on our experience, the averaged weights present good repeatability, whereas in each individual training different weights may be achieved (even if higher number of data is used for an individual training). To make three controllers comparable, they are tuned so that they have similar rise time. Also, experiments are done with $500 Hz$ sampling frequency. Furthermore, no friction compensation is done in experiments. For all experiments, actuators are saturated at $5 \, N.m$. Other parameters used in the tests can be found in~\cref{appendixparams}.

\subsection{Results}
Two scenarios are considered to compare the performance of the methods. First, the robot is supposed to draw a circle in $z$-plane, with $x = 250\, cos(\pi t)\, (mm)$ and $y = 250\, sin(\pi t)\, (mm) $. Second, the robot is supposed to go to two different locations sequentially, i.e., $\begin{bmatrix}100& 100& 450\end{bmatrix}^\intercal \,(mm)$ and \linebreak $\begin{bmatrix}-100& -100& 600\end{bmatrix}^\intercal \,(mm)$. Moreover, to compare robustness of the controllers, both scenarios are repeated by adding a $1kg$ mass as an uncertainty to the end effector. Also all experiments started from robot's home position at $\begin{bmatrix}0& 0& 500\end{bmatrix}^\intercal \,(mm)$. Video of the tests can be found in \cite{testvideo}.

\begin{figure}[H]
	\centering
	\includegraphics[width=.4\textwidth]{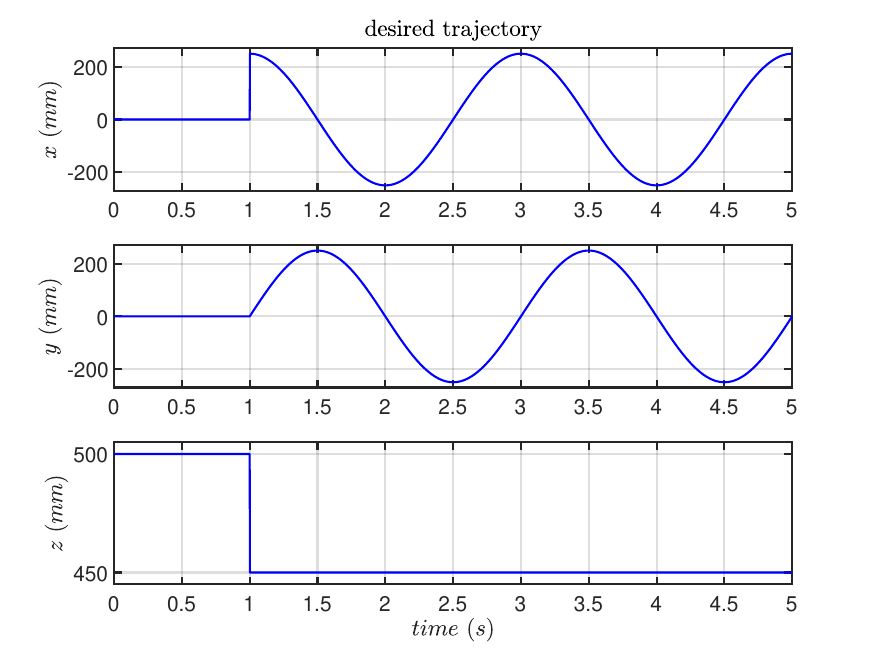}
	\caption{First scenario desired trajectory}\label{fig:circtraj}
\end{figure}

\begin{figure}[H]
	\centering
	\includegraphics[width=.4\textwidth]{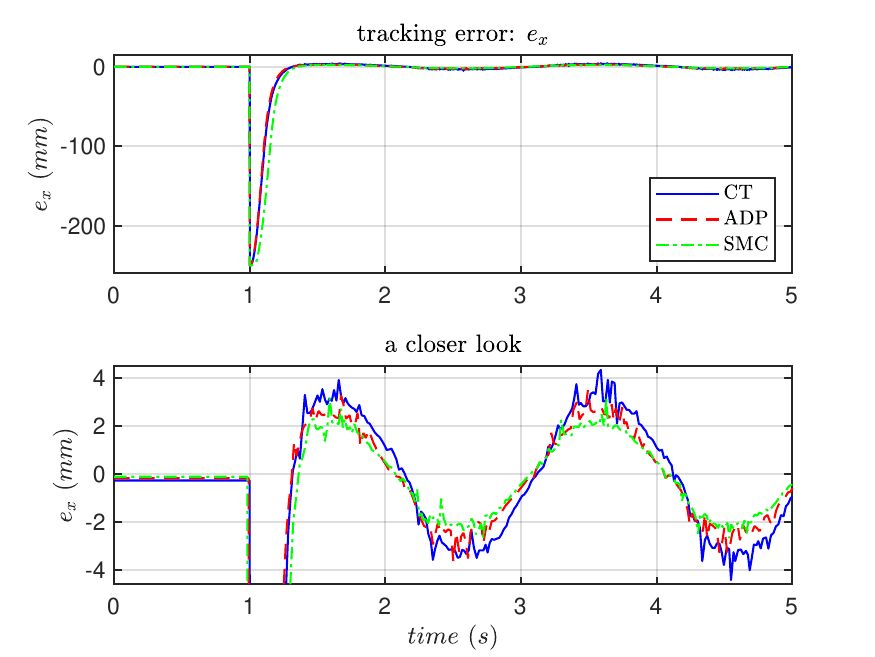}
	\caption{First scenario tracking error $e_x$: without uncertainty}\label{fig:excirc}
\end{figure}

\begin{figure}[!htp]
	\centering
	\includegraphics[width=.48\textwidth]{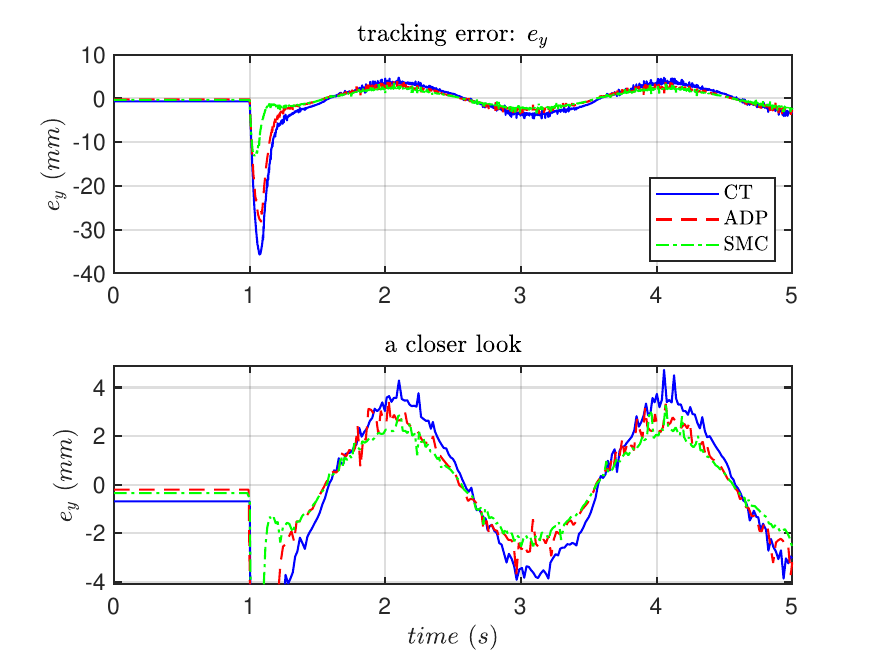}
	\caption{First scenario tracking error $e_y$: without uncertainty}\label{fig:eycirc}
\end{figure}

\begin{figure}[!htp]
	\centering
	\includegraphics[width=.48\textwidth]{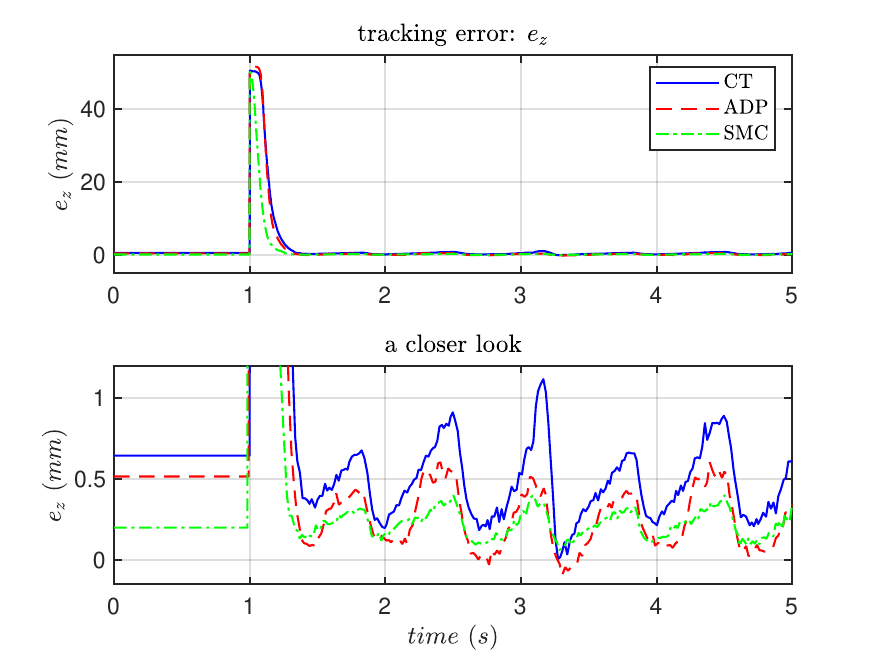}
	\caption{First scenario tracking error $e_z$: without uncertainty}\label{fig:ezcirc}
\end{figure}

Results related to first scenario without uncertainty are summarized in \cref{fig:circtraj,fig:excirc,fig:eycirc,fig:ezcirc,table:ms}.
\begin{table}[H]
	\centering
	\begin{tabular}{ l c c c  }
		{} & x& y& z\\
		\hline
		$\mu_{\left|e\right|}\,(mm)$ & & &\\ \hline\hline
		CT & $2.2119$ & $2.1916$ & $0.4873$\\
		ADP & $1.6669$ & $1.7063$ & $0.2387$\\
		SMC & $1.4699$ & $1.4020$ & $0.2214$\\\hline
		$\sigma_{\left|e\right|}\,(mm)$ & & &\\ \hline\hline
		CT & $1.1095$ & $1.1505$ & $0.2358$\\
		ADP & $0.8907$ & $0.8686$ & $0.1751$\\
		SMC & $0.7393$ & $0.7247$ & $0.0861$
	\end{tabular}
	\caption{Mean and standard deviation of steady state $\left|e\right|$,\\first scenario without uncertainty}\label{table:ms}
\end{table}

\begin{figure}[H]
	\centering
	\includegraphics[width=.48\textwidth]{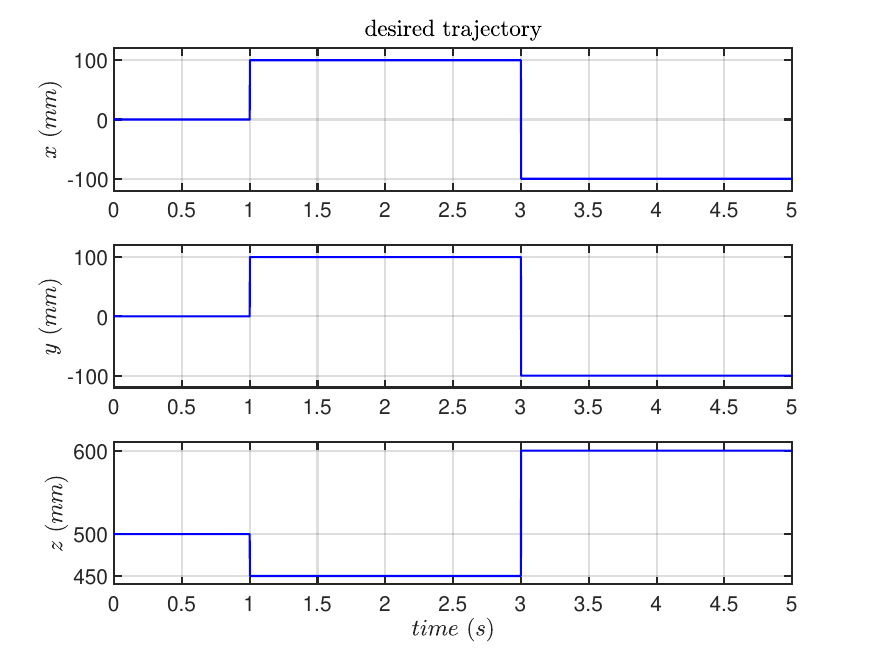}
	\caption{Second scenario desired trajectory}\label{fig:steptraj}
\end{figure}

\begin{figure}[H]
	\centering
	\includegraphics[width=.48\textwidth]{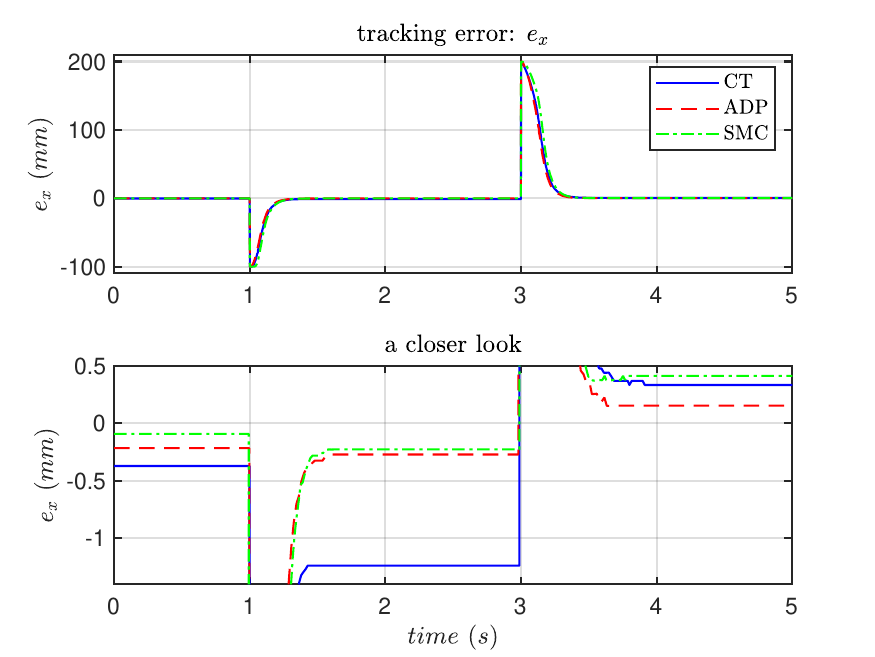}
	\caption{Second scenario tracking error $e_x$: without uncertainty}\label{fig:exstep}
\end{figure}

\begin{figure}[H]
	\centering
	\includegraphics[width=.48\textwidth]{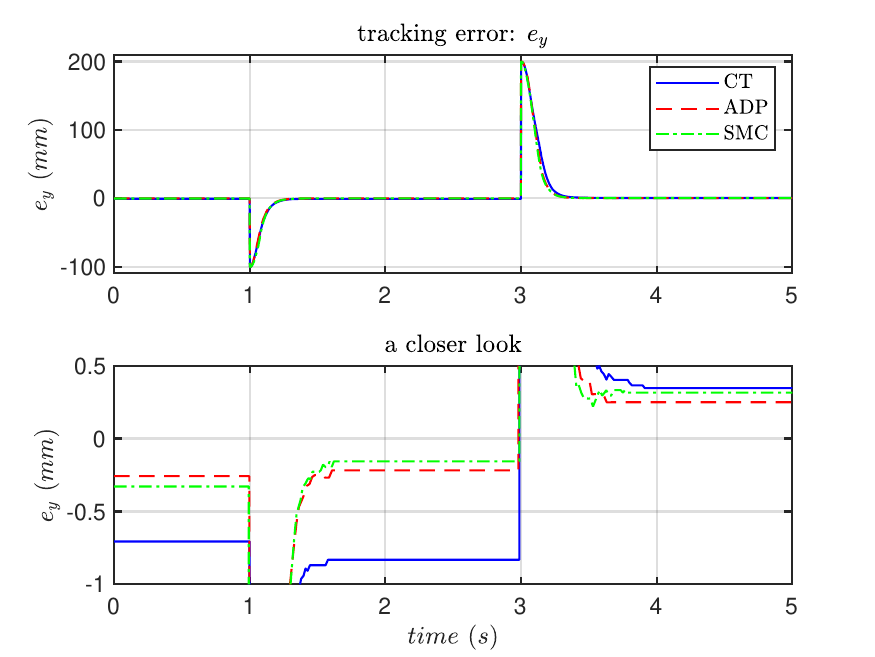}
	\caption{Second scenario tracking error $e_y$: without uncertainty}\label{fig:eystep}
\end{figure}

\begin{figure}[H]
	\centering
	\includegraphics[width=.48\textwidth]{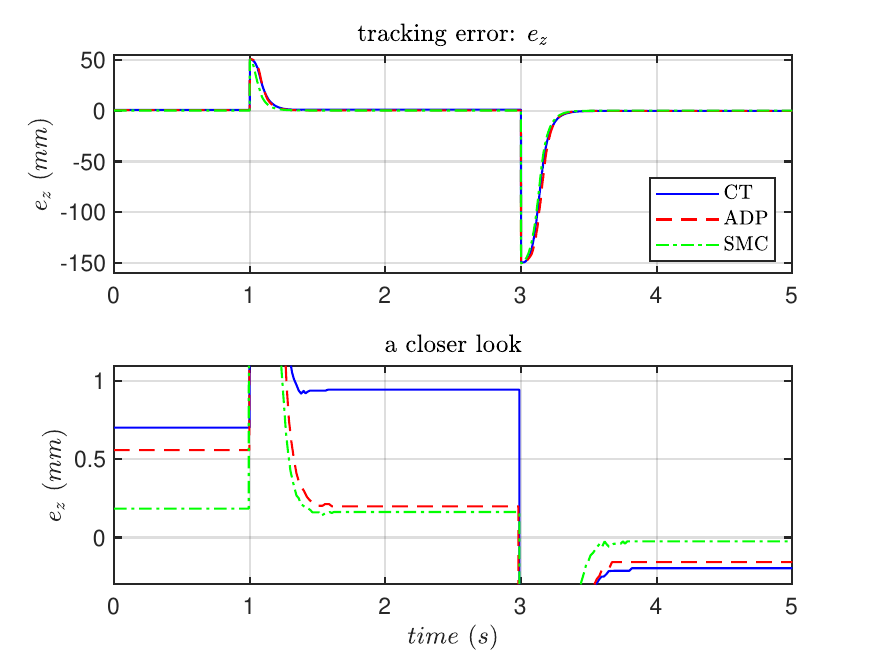}
	\caption{Second scenario tracking error $e_z$: without uncertainty}\label{fig:ezstep}
\end{figure}

Results related to second scenario without uncertainty are presented in \cref{fig:steptraj,fig:exstep,fig:eystep,fig:ezstep}. The difference between performance of the three methods (without the uncertain mass), as observed through \cref{fig:excirc,fig:eycirc,fig:ezcirc,fig:exstep,fig:eystep,fig:ezstep,table:ms}, is small. However, even these small differences are considerable in the context of robotic applications, given tight tolerances and high accuracy requirements. CT does the worst, both in step response and following a circle. The performance of the proposed method is very close to that of SMC. But except for the $y$ coordinate error of step test, SMC controller slightly does a better job.

\begin{figure}[H]
	\centering
	\includegraphics[width=.48\textwidth]{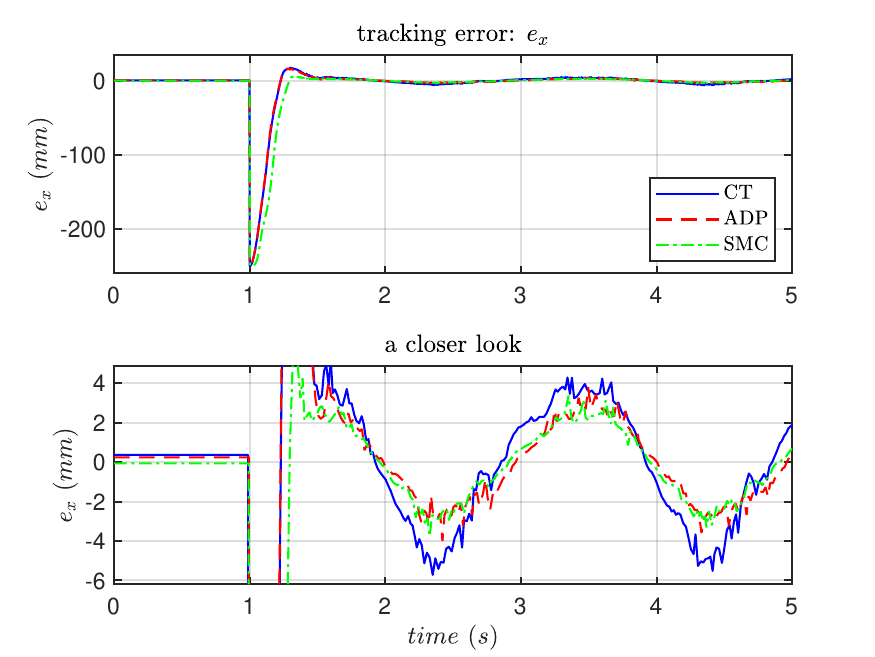}
	\caption{First scenario tracking error $e_x$: with uncertainty}\label{fig:excircu}
\end{figure}

\begin{figure}[H]
	\centering
	\includegraphics[width=.48\textwidth]{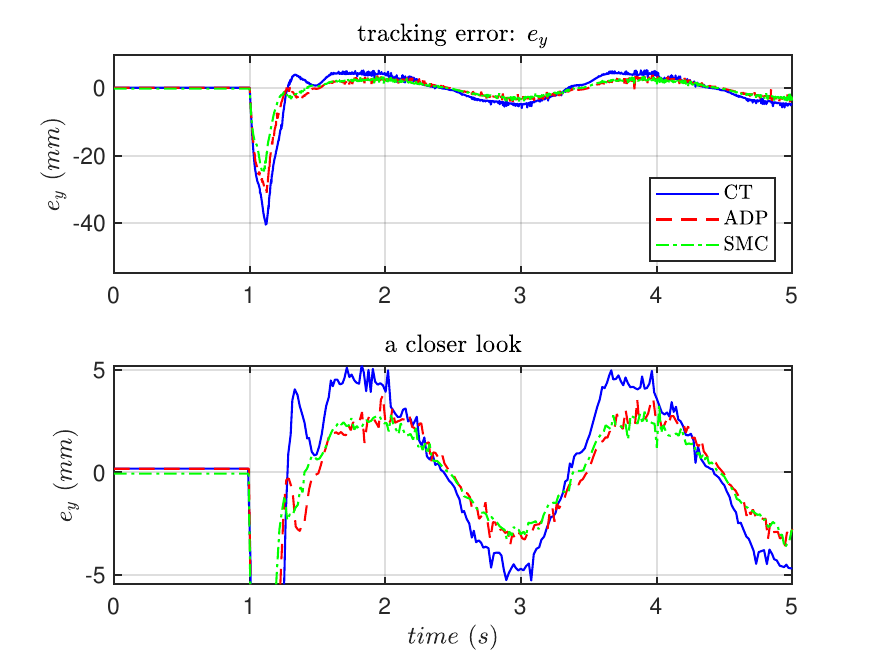}
	\caption{First scenario tracking error $e_y$: with uncertainty}\label{fig:eycircu}
\end{figure}

\begin{figure}[H]
	\centering
	\includegraphics[width=.48\textwidth]{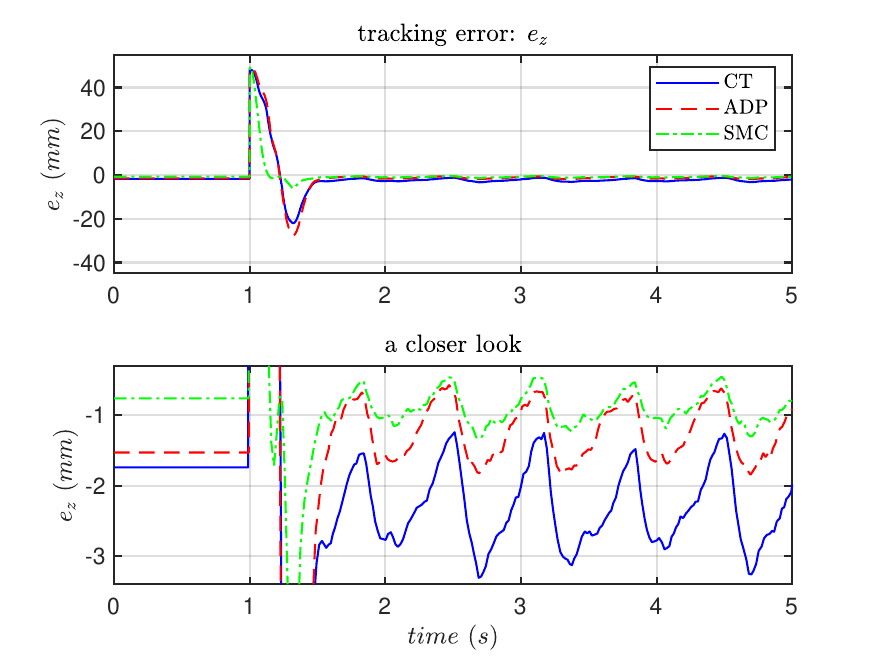}
	\caption{First scenario tracking error $e_z$: with uncertainty}\label{fig:ezcircu}
\end{figure}

\begin{figure}[H]
	\centering
	\includegraphics[width=.48\textwidth]{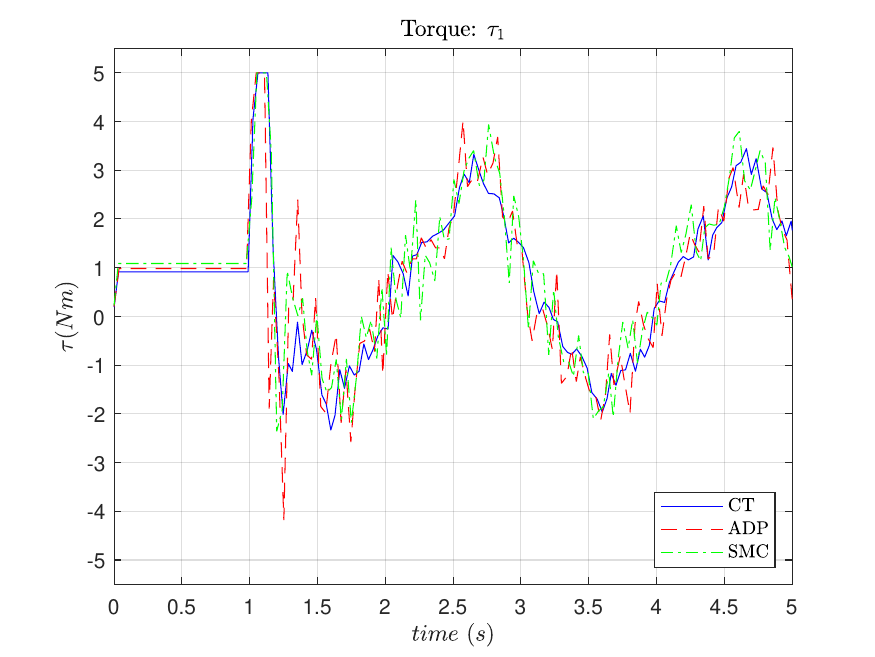}
	\caption{First scenario torque $\tau_1$: with uncertainty}\label{fig:tauxcircu}
\end{figure}

\begin{figure}[!tp]
	\centering
	\includegraphics[width=.48\textwidth]{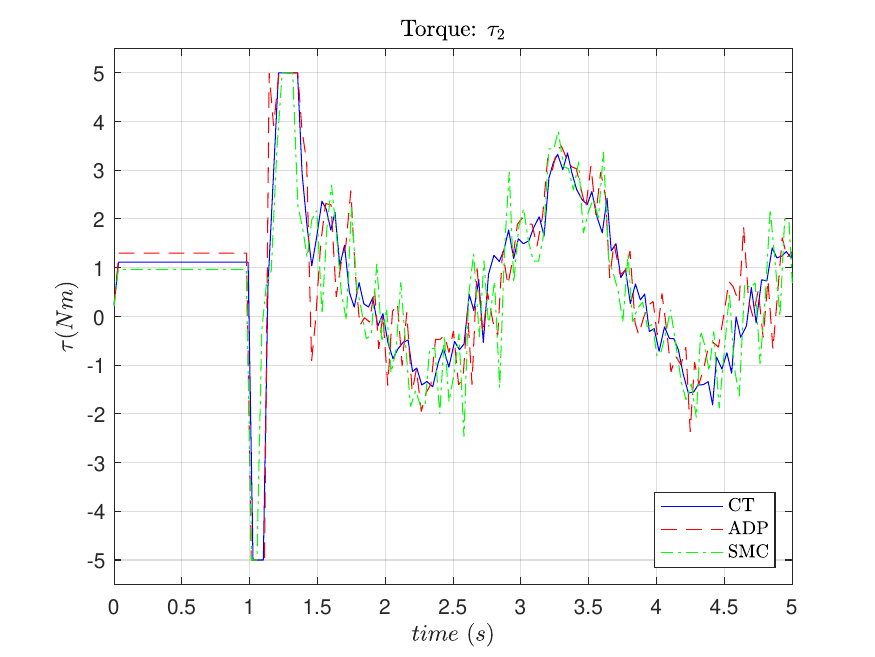}
	\caption{First scenario torque $\tau_2$: with uncertainty}\label{fig:tauycircu}
\end{figure}

\begin{figure}[!tp]
	\centering
	\includegraphics[width=.48\textwidth]{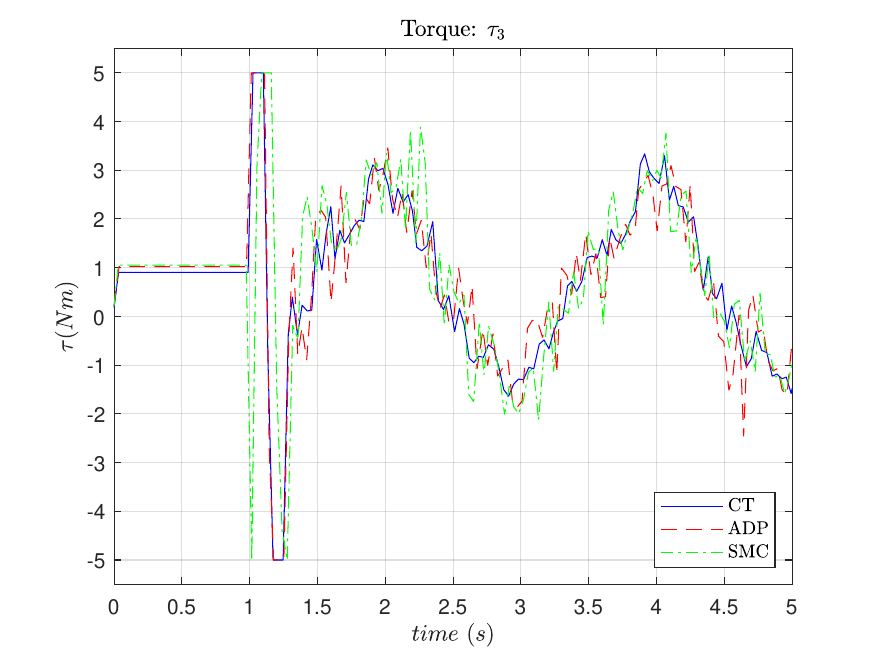}
	\caption{First scenario torque $\tau_3$: with uncertainty}\label{fig:tauzcircu}
\end{figure}

\begin{figure}[!tp]
	\centering
	\includegraphics[width=.48\textwidth]{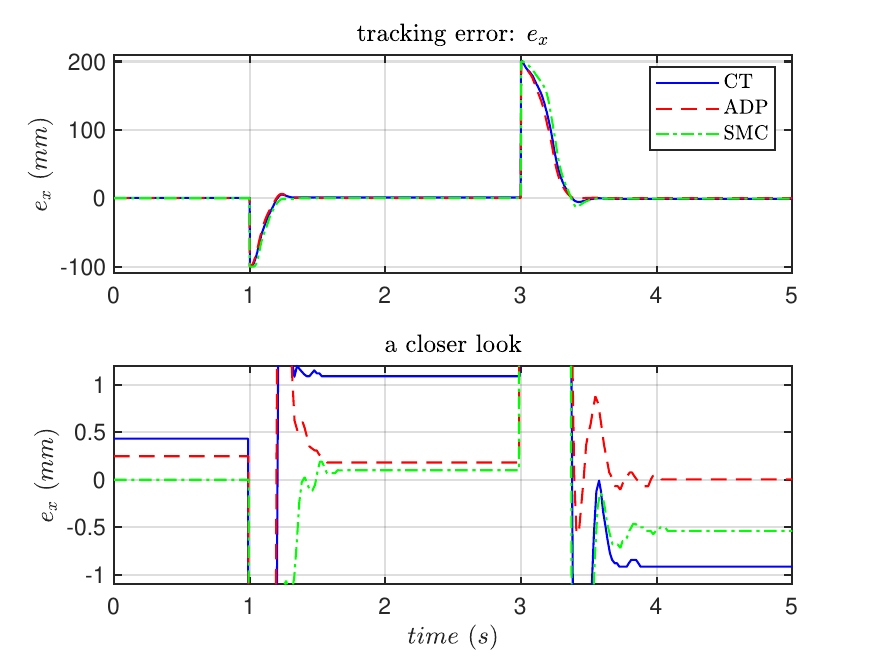}
	\caption{Second scenario tracking error $e_x$: with uncertainty}\label{fig:exstepu}
\end{figure}

\begin{figure}[!tp]
	\centering
	\includegraphics[width=.48\textwidth]{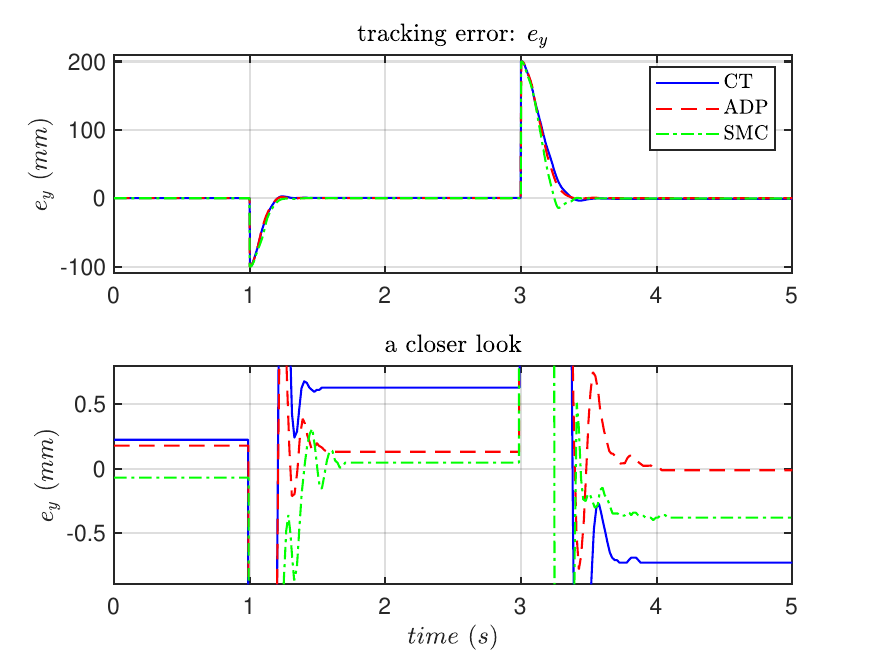}
	\caption{Second scenario tracking error $e_y$: with uncertainty}\label{fig:eystepu}
\end{figure}

\begin{figure}[!tp]
	\centering
	\includegraphics[width=.48\textwidth]{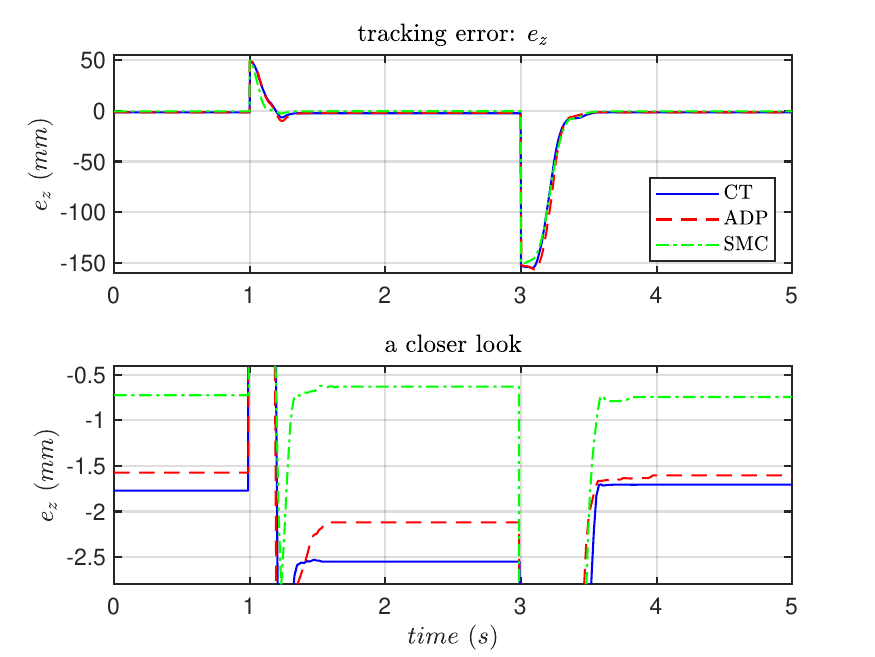}
	\caption{Second scenario tracking error $e_z$: with uncertainty}\label{fig:ezstepu}
\end{figure}

\begin{figure}[!tp]
	\centering
	\includegraphics[width=.48\textwidth]{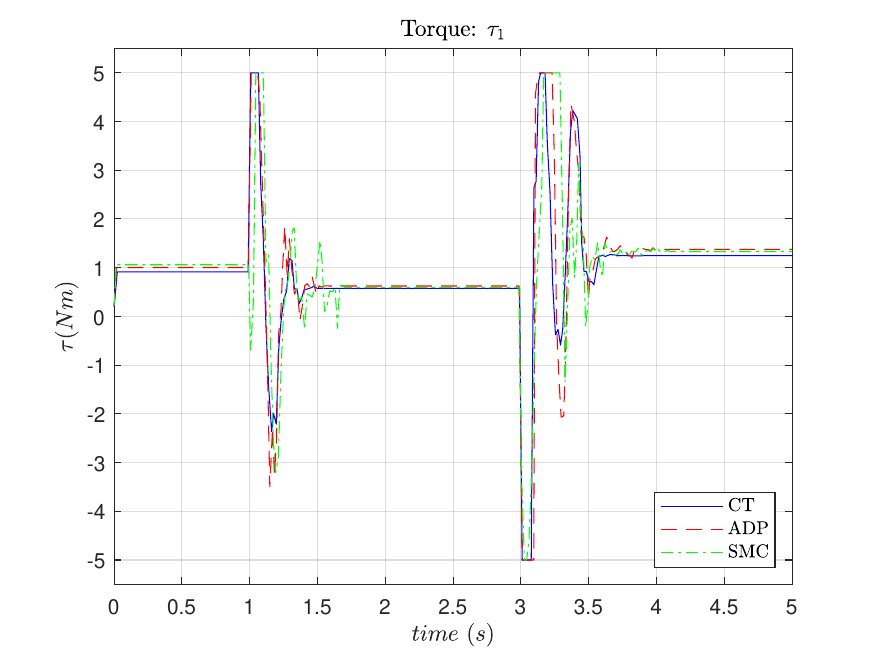}
	\caption{Second scenario torque $\tau_1$: with uncertainty}\label{fig:tauxstepu}
\end{figure}

\begin{figure}[!tp]
	\centering
	\includegraphics[width=.48\textwidth]{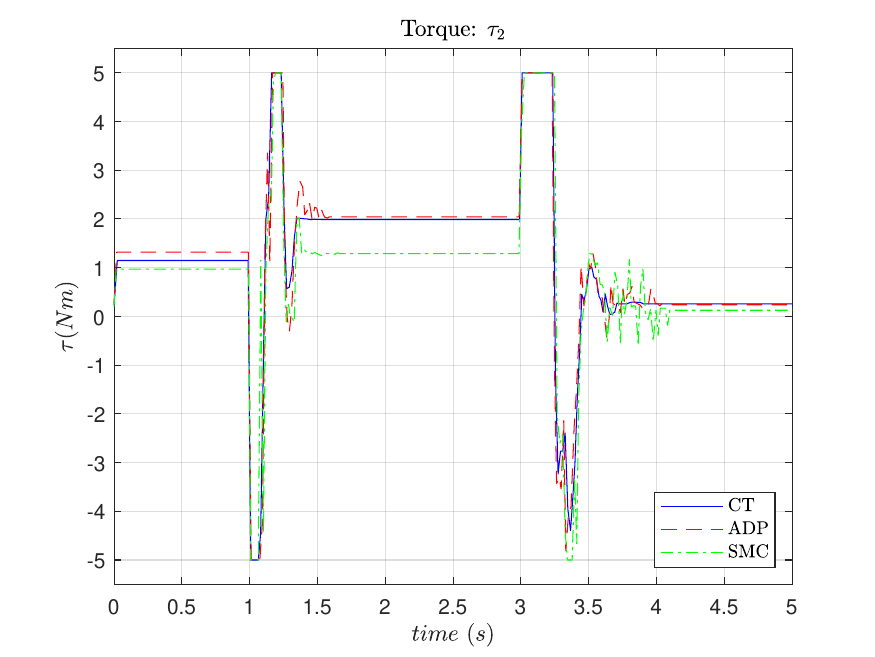}
	\caption{Second scenario torque $\tau_2$: with uncertainty}\label{fig:tauystepu}
\end{figure}

\begin{figure}[!tp]
	\centering
	\includegraphics[width=.48\textwidth]{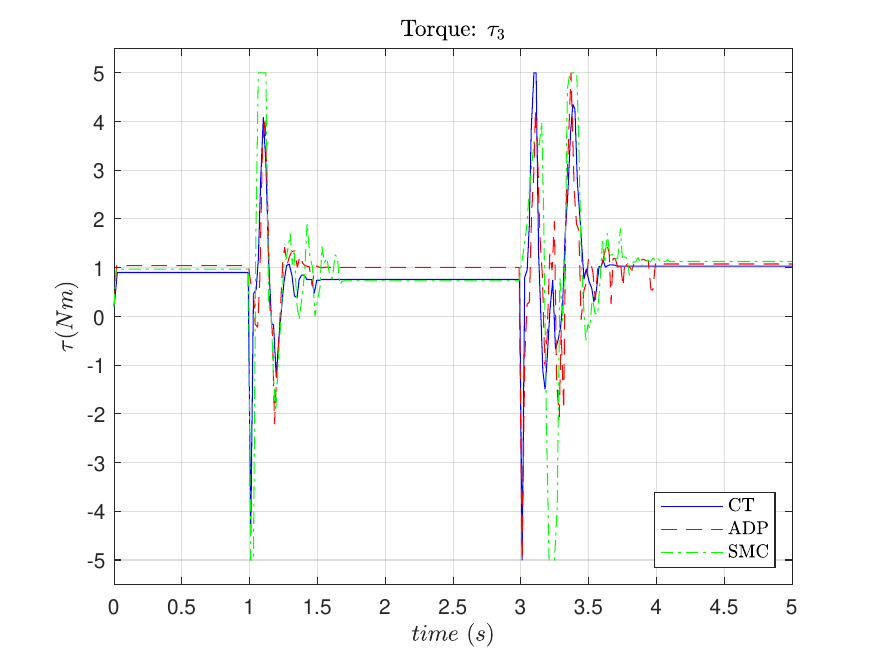}
	\caption{Second scenario torque $\tau_3$: with uncertainty}\label{fig:tauzstepu}
\end{figure}

Results of experiments with a mass added as uncertainty are shown
\cref{fig:excircu,fig:eycircu,fig:ezcircu,fig:tauxcircu,fig:tauycircu,fig:tauzcircu,table:msu,fig:exstepu,fig:eystepu,fig:ezstepu,fig:tauxstepu,fig:tauystepu,fig:tauzstepu}.
\begin{table}[H]
	\centering
	\begin{tabular}{ l c c c  }
		{} & x& y& z\\
		\hline
		$\mu_{\left|e\right|}\,(mm)$ & & &\\ \hline\hline
		CT & $2.6833$ & $2.8572$ & $2.2906$\\
		ADP & $1.7528$ & $1.8637$ & $1.2390$\\
		SMC & $1.6172$ & $1.7284$ & $0.8945$\\\hline
		$\sigma_{\left|e\right|}\,(mm)$ & & &\\ \hline\hline
		CT & $1.4124$ & $1.5423$ & $0.55518$\\
		ADP & $0.9451$ & $0.9336$ & $0.3820$\\
		SMC & $0.8496$ & $0.8882$ & $0.2250$
	\end{tabular}
	\caption{Mean and standard deviation of steady state $\left|e\right|$,\\first scenario with uncertainty}\label{table:msu}
\end{table}


As it can be seen in \cref{fig:excircu,fig:eycircu,fig:ezcircu,fig:exstepu,fig:eystepu,fig:ezstepu,table:msu}, adding the uncertainty to the system increases steady state errors for all three methods. While one can say the robustness of SMC is higher than other two methods, the performance of the proposed methods still remains close to that of SMC. Moreover, the computed torque controller falls much behind in comparison to the other two methods, as expected.

The cost comparison of the controllers is given for first scenario, with and without uncertainty in \cref{table:cost}. This cost is evaluated based on considering the total cost in the objective function, similar to \cref{eq:ctotcost}, for first five seconds of the experiments. As expected, for the case without uncertainty, the cost of the proposed approach is lower by at least $22\%$. In the case with uncertainty, the cost of proposed method is approximately $2\%$ higher than that of computed torque method. This is logical as the superior performance of the proposed method against uncertainty, needs more control effort.
\begin{table}[!tp]
	\centering
	\begin{tabular}{ l c c c  }
		{} & CT& ADP& SMC\\
		\hline\hline
		Without uncertainty& $681.65$&$558.43$&$1180.49$\\ \hline
		With uncertainty& $828.62$&$851.45$&$1181.44$\\\hline\hline
	\end{tabular}
	\caption{Total cost for first scenario}\label{table:cost}
\end{table}

Despite advantages of the proposed optimal controller, that can be seen in experiments, like every other method it has its disadvantageous too. ADP, that is used to solve the nonlinear optimal control problem in a closed form, is a numerical procedure and it can have convergence problems in practice. To be more precise, depending on the approximation error of the neural network that is used, the ADP based algorithm may not converge for all values of $Q,R$, and sampling time. Therefore, the designer should consider this in controller design stage. The other disadvantage of ADP is about choosing basis functions. Even though some works has been done about this, there is no conclusive work yet. Generally as ROI gets bigger in dimension, finding basis functions that accurately interpolate the value function, gets much harder. Consequently the convergence of the algorithm will be affected, however normalizing the data might be helpful. The other issue is the so called curse of dimensionality \cite{Bellman:1957} in dynamic programming. Despite mitigation of this problem by implementing ADP, and also in the proposed method through reduction of value function parameters by introducing expectation of value function (expectation of value function depends only on $E$ instead of $E$ and $X_d$), the problem still exists.

The other point to be mentioned is about chattering. This is usually considered a problem related to SMC. However, this phenomenon can also happen for the other two methods, because of control discontinuity resulted from digital implementation. In tuning all of the controllers, it was observed that there is an upper-bound on the aggressiveness of each of the controllers, that can be achieved without chattering.
\section{Conclusion}
In this paper, a new framework is introduced for optimal tracking problem of a class of nonlinear systems. In contrast to previous works on optimal control, the presented approach can track any trajectories (of course in the ROI) after one training. Also using expectation of total cost, number of parameters decreased. This mitigates curse of dimensionality. The presented method is then applied to a relatively complex nonlinear system and its performance is shown experimentally. The current work addressed asymptotic optimal tracking problem for nonlinear systems in canonical form. Future work can focus on the same problem for general input affine systems with well defined relative degree.
\section{Acknowledgments}
This research was partially supported by the United States National Science Foundation through Grant 1745212.
\crefalias{section}{appendix}
\section*{Appendix}
\setcounter{section}{0}
\renewcommand*{\theHsection}{Appendix.\the\value{section}}
\renewcommand\thesubsection{\Alph{subsection}}
\crefalias{subsection}{appendix}
\subsection{}\label{appendixexample}
Consider the following optimal tracking system that is defined based on exact total cost
\begin{equation*}
\dot{x} = x + u,
\end{equation*}
\begin{equation*}
J(x_0) = \int_{0}^{\infty} exp(-\rho t)(q(x-r)^2 + u^2) dt,
\end{equation*}
where $\rho>0$ and $q = 1$. Let assume that the desired trajectory is $r(t) = 2$. In this case if the initial condition is $x_0 = 1$. Then the optimal control intuitively becomes $u^*(t) = -1$, this can be easily verified from the standard LQT solution. Applying the optimal control, the state time history becomes $x(t) = 1$.  Therefore, the optimal tracking control based on exact total cost is not asymptotically stabilizing. This happens due to the fact that the reference trajectory is not within invariant sets of the system (this makes error dynamics non-stationary at origin), which is needed for asymptotic stability of the closed loop system under optimal control. To be more precise the optimal tracking controller based on the total cost can only asymptotically track reference trajectories that are among invariant sets of the system. This solution is acceptable in an economical optimization problem, however in control context, asymptotic stability is favored. To achieve general asymptotic stability from optimal control resulted from optimizing exact total cost one typically needs $q \rightarrow \infty$.
\subsection{Experimental parameters}\label{appendixparams}
Here, control parameters and geometrical and inertial parameters of the robot used in the experiments are given.
\noindent \newline
Computed torque parameters:
\begin{equation*}
\begin{array}{lr}
K_p = 1600 & \text{proportional gain}\\
K_d = 100 & \text{derivative gain}
\end{array}
\end{equation*}

\noindent \newline
Proposed controller parameters:
\begin{equation*}
\begin{array}{ll}
Q = D^{\intercal} D & \text{state penalizing matrix}
\end{array}
\end{equation*}
\begin{equation*}
\begin{array}{ll}
\text{where} & D = \begin{bmatrix}
20 & 0 & 0 & 1 & 0 & 0\\
0 & 20 & 0 & 0 & 1 & 0\\
0 & 0 & 20 & 0 & 0 & 1
\end{bmatrix}
\end{array}
\end{equation*}
\begin{equation*}
\begin{array}{ll}
R = 0.001I_{3 \times 3} & \text{control penalizing matrix}\\ 
\end{array}
\end{equation*}
\begin{equation*}
\begin{array}{ll} 
\begin{aligned}
\varphi(X) = \left[x_1x_2 \quad x_1x_3  \quad x_1x_4  \quad \dotso \right.\\
\qquad x_2x_3  \quad x_1x_5  \quad x_2x_4 \\
\qquad x_1x_6  \quad x_2x_5  \quad x_3x_4 \\
\qquad x_2x_6  \quad x_3x_5  \quad x_3x_6 \\
\qquad x_4x_5  \quad x_4x_6  \quad x_5x_6 \\
\qquad {x_1}^2  \quad {x_2}^2  \quad {x_3}^2 \\
\left. \qquad {x_4}^2  \quad {x_5}^2  \quad {x_6}^2 \right]^{\intercal}\end{aligned}
& \text{basis function}
\end{array}
\end{equation*}
where
\begin{multline*}
\begin{bmatrix}
x_1 \quad x_2 \quad x_3
\end{bmatrix} = 
\begin{bmatrix}
e_x \quad e_y \quad e_z
\end{bmatrix} \text{\qquad and}\\ \begin{bmatrix}
x_4 \quad x_5 \quad x_6
\end{bmatrix} = 
\begin{bmatrix}
\dot{e}_x \quad \dot{e}_y \quad \dot{e}_z
\end{bmatrix}
\end{multline*}
\begin{multline*} 
\begin{aligned}
W = \left[0.0025 \quad -0.1939  \qquad 0.0330  \qquad \dotso \right.\\
\quad -0.2257  \qquad 0.0026  \quad -0.0009 \\
\qquad 0.0008  \qquad 0.0317  \quad -0.0026 \\
\qquad 0.0002  \quad -0.0055  \qquad 0.0507 \\
\quad -0.0001  \quad -0.0002  \quad -0.0002 \\
\qquad 1.8550  \qquad 1.8911  \qquad 1.9928 \\
\left. \qquad 0.0012  \qquad 0.0012  \qquad 0.0016 \right]^{\intercal}\end{aligned}
\\ \text{optimal weights}
\end{multline*}

\noindent \newline
Sliding mode parameters:
\begin{equation*}
\begin{array}{ll}
K = 70 & \text{sliding mode gain}\\
\lambda = 20 & \text{sliding surface parameter}\\
\varphi = 0.35 & \text{boundary layer}
\end{array}
\end{equation*}
Geometrical parameters:
\begin{equation*}
\begin{array}{ll}
r_b = 0.2 \: m & \text{fixed platform radius}\\
r_a = 0.05 \: m & \text{moving platform radius}\\
\varphi_1 = \pi & \text{1st leg position angle}\\
\varphi_2 = -\frac{\pi}{3} & \text{2nd leg position angle}\\
\varphi_3 = \frac{\pi}{3} & \text{3rd leg position angle}\\
l_L = 0.2 \: m & \text{motor leg length}\\
r_L = 0.3933 & \text{motor leg COM lenght ratio}\\
l_K  = .52 \: m & \text{intermediary leg length}\\
r_K = .5 & \text{intermediary leg COM length ratio}\\
R_b = \begin{bmatrix}
0 & r_b & 0
\end{bmatrix}^\intercal\\
R_a = \begin{bmatrix}
0 & r_a & 0
\end{bmatrix}^\intercal
\end{array}
.
\end{equation*}
Inertial parameters: \newline
\begin{tabular}{ll}
	$m_P = 1.055 \: kg$ & moving platform mass\\
	$I_{mo} = 0.0465475 \: kg.m^2$ & motor inertia\\
	$m_L = 0.116 \: kg$ & motor leg mass\\
	$m_K = 2\times 0.05788 \: kg$ & intermediary leg mass\\
	$I_L = 6.4345319\times 10^{-4} \: kg.m^2$ & motor leg inertia\\
	\multirow{ 3}{*}{$d_K = 0.08 \: m$} & intermediary leg distance\\
	&from each other at every\\
	& joint\\
	\multirow{ 3}{*}{$I_K = 5.74769459 \times 10^{-3} \: kg.m^2$} & intermediary leg pairs\\
	&dominant inertia\\
\end{tabular}
\bibliographystyle{imaiai}

\end{document}